\documentclass[reqno]{amsart}
\usepackage{amssymb, amsmath, amsthm}

\def\zz{{\bf Z}}
\def\qq{{\bf Q}}

\def\cs{\mathop{\#}}
\newcommand\rk{\operatorname{rk}}

\newtheorem{theorem}{Theorem}
\newtheorem{lemma}[theorem]{Lemma}

\theoremstyle{definition}
\newtheorem{definition}[theorem]{Definition}

\def\calr{\mathcal{R}}
\def\cals{\mathcal{S}}

\numberwithin{equation}{section}

\def\co{\colon}

\begin{document}

\title{ The Hausmann-Weinberger 4--Manifold Invariant of Abelian Groups}

\author{Paul Kirk}

\thanks{The first named author gratefully acknowledges the support of
the
National Science Foundation under   grant no. DMS-0202148.}

\author{Charles Livingston}

\address{Department of Mathematics, Indiana University, Bloomington, IN 47405}

\email{pkirk@indiana.edu}
\email{livingst@indiana.edu}

\keywords{ Hausmann-Weinberger invariant, fundamental group,
four-manifold, minimal Euler characteristic}

\subjclass{}

\date{\today}

\begin{abstract}  The Hausmann-Weinberger invariant of a group $G$ is the minimal
Euler characteristic of a closed orientable 4--manifold $M$ with fundamental
group $G$.  We compute this invariant for finitely generated free abelian groups
and estimate the invariant  for all finitely generated abelian groups.
\end{abstract}

\maketitle

\vskip.2in
\section{Introduction}

For any finitely   presented group $G$ there exists a closed oriented 4--manifold
$M$ with $\pi_1(M) = G$.   Hausmann and Weinberger
\cite{hw} defined the integer valued invariant
$q(G)$ to be the least Euler characteristic among all such $M$.  The explicit
construction of a 4--manifold with   $\pi_1(M) = G$, based on a presentation of $G$,
 yields an
upper bound on $q(G)$.  As  pointed out in~\cite{hw}, the isomorphism $H_1(M)\to
H_1(G)$, the surjection $ H_2(M)\to H_2(G)$, and  Poincar\'e duality yield  a
lower bound.  Together these bounds are
\begin{equation}\label{lbound} 2-2\beta_1(G)+\beta_2(G)\leq q(G) \le 2-
2\mbox{def}(G),
\end{equation} where $\mbox{def}(G)$ is the deficiency of $G$, the maximum
possible difference $g-r$ where the $g$ is the number of generators and
$r$ the number of relations in a presentation of
$G$, and $\beta_i(G)$ denotes the $i$th Betti number of $G$ (with some fixed
coefficients).

Since~\cite{hw}, advances have been made in the study of this invariant, most
notably through the methods of  $l^2$--homology.  For instance, in~\cite{e1, e2}
Eckmann proves   that for infinite amenable groups $G$,
$q(G) \ge 0$.  L\"{u}ck~\cite{luck} extended this to all group $G$ with
$b_1(G) =0$, where $b_1$ denotes the first $l^2$--Betti number.  Other work
includes~\cite{jk} and especially the paper by Kotschick~\cite{ko} in which
Problem 5.2 asks for the explicit value of $q(\zz^n)$. The general problem
 of computing $q(G)$ appears as Problem 4.59 in
Kirby's problem list,~\cite{ki}.

Despite these past efforts, the  Hausmann-Weinberger invariant remains uncomputed
for some of the most elementary groups.  In ~\cite{hw} it is  observed that
$q(\zz^n)$ is given by 2, 0, 0, 2, and 0, for $n= 0, 1, 2, 3, $ and $4$,
respectively. (For the case of $n=3$,~\cite{hw} refers to an unpublished argument
of  Kreck.  Proofs appear  in~\cite{e1,  ko}.)   If $\Gamma_g$ denotes the
fundamental group of a surface of genus $g$,~\cite{ko} computes $q(\Gamma_g)$
and  $q(\Gamma_{g_1} \times \Gamma_{g_2})$.  If a closed 4-manifold $X$ is aspherical then
$\chi(X)=q(\pi_1(X))$.
  Beyond this, few explicit values
of $q(G)$ have been calculated.  Our main theorem is the following:

\begin{theorem}\label{mainthm}  With the exceptions of $q(\zz^3) = 2$ and
$q(\zz^5) = 6$,
$q(\zz^n)$ is given by:

\begin{displaymath} q(\zz^n) = \begin{cases} (n-1)(n-4)/2, & \mbox{if\ \ } n
\equiv 0 \mbox{\ or\ } n \equiv 1 \mod 4; \\ (n-2)(n-3)/2, & \mbox{if\ \ } n
\equiv 2
\mbox{\ or\ } n \equiv 3 \mod 4. \\
\end{cases}
\end{displaymath}
\end{theorem}

For contrast, the bounds~(\ref{lbound}) give only that
$$(n-1)(n-4)/2 \le q(\zz^n) \le (n-1)(n-2).$$ It is straightforward to check that
an alternative way to state Theorem~\ref{mainthm} is that for
$n\ne  3,5$, the lower bound of (\ref{lbound}) is attained  when
the binomial coefficient $C(n,2)$ is even, and if $C(n,2)$ is odd the value is one more than the lower
bound of  ~(\ref{lbound}).

The calculations of $q(\zz^n)$ can be used to estimate (and in some cases
calculate) $q(G)$ for  other groups. We examine the problem for finitely
generated abelian  groups and prove the following. In the statement, $\epsilon_n$
equals zero or one according to whether $C(n,2)$ is even or odd, respectively.
\vskip.05in

\begin{theorem} \label{genral} Let $G=\zz/d_1\oplus \cdots\oplus
\zz/d_k\oplus \zz^n$ with $d_i|d_{i+1}$. Suppose that  $k\ge 1$ and   $ k + n  \ne
3,4,5$ or $6$. Then
  $$  0\leq q(G)-\big( 1-n+C(n+k-1,2)\big)\leq
\operatorname{min}\{|n-1|+\epsilon_{n+k-1}, k+\epsilon_{n+k}    \}.$$
  Moreover, $q(\zz/d)=2$,  $q(\zz/d_1\oplus\zz/d_2)=2$,  and  if   $C(k,2)$ is even
and $k\ne 5$,  then
  $$q(\zz/d_1\oplus \cdots \oplus \zz/d_k\oplus \zz)= C(k,2).$$
  \end{theorem}

In closing this introduction we   mention results concerning the evaluation
of $q(P)$ where $P$ is a perfect group.  Here the bounds given by~(\ref{lbound}) are $$2
+
\beta_2(P)
\le q(P)
\le 2 -
2\mbox{def}(G).$$  In~\cite{hw} perfect groups $P$ are
constructed with $\beta_2(P) = 0$ but $q(P) > 2$.   Hillman~\cite{hi}
constructed perfect groups of deficiency
$-1$ with $q(P) = 2$ and the second author~\cite{li}   extended this  
  to find perfect $P$ with arbitrarily large negative deficiency and $q(P) =
2$.

\section{Notation and basic results} A slightly different invariant,
$h(G)$, can be defined to  be the minimum value of
$\beta_2(M)$ among all  oriented closed 4-manifolds $M$ with
$\pi_1(M)=G$.  We abbreviate
$h(\zz^n) = h(n)$. Clearly
$q(G) = 2 - 2\beta_1(G) + h(G)$, so the invariants are basically equivalent.  It
is  more convenient here to work in terms of $h$.  The bounds (\ref{lbound}) on
$q(\zz^n)$ translate to the bounds
$$C(n,2)\leq h(n)\leq 2C(n,2).$$
We introduce the following   auxiliary function:

\begin{displaymath}
\epsilon_n = \begin{cases} 0,  & \mbox{if\ \ } C(n,2) \mbox{\ is even;}
\\ 1, &
\mbox{if\ \ }  C(n,2) \mbox{\ is odd}.
\end{cases}
\end{displaymath}
Since $C(n,2)$ is even if and only if $n \equiv 0$ or $1$ mod 4, 
Theorem~\ref{mainthm} can be restated as follows.

\vskip.1in
\noindent{\bf Theorem~\ref{mainthm}.} {\em With the exceptions of $h(3) = 6$ and
$h(5) = 14$, $h(n) = C(n,2)+\epsilon_n$ for all $n$.}
\vskip.05in

Basic examples of 4--manifolds will be built from products of surfaces.  For $n$
even we will denote by $F_n$ the closed orientable surface of genus $n/2$.


\section{Bounds on $h(n)$}

\begin{theorem}  If $h(n) = C(n,2)$, then $C(n,2)$ must be even.  Thus
$h(n) \ge C(n,2)+\epsilon_n$.
\end{theorem}

\begin{proof}    If $\phi\co \pi_1(M) \to \zz^n$, we have $\phi_* \co H_*(M) \to
H_*(\zz^n)$.  Dually there is the map of cohomology rings
$\phi^* \co H^*(\zz^n) \to H^*(M)$.  Notice that $H^*(\zz^n)$ is an exterior algebra
on the generators $e_1,\cdots, e_n\in H^1(\zz^n)$.

Suppose $\phi\co \pi_1(M) \to \zz^n$ is an isomorphism and
$\beta_2(M) = C(n,2)$. Then the map $\phi_2 \co H_2(M) \to  H_2(\zz^n)$ is a
surjection from $\zz^{C(n,2)}$ to $\zz^{C(n,2)}$, and hence is an isomorphism.
It follows that $\phi^2 \co H^2(\zz^n)
\to H^2(M)$ is also   an isomorphism.

  Since
$(e_i  e_j)^2=0$, $H^2(M)$ has a basis for which all squares are
zero.  It follows that the intersection form of
$M$ is even.  But even unimodular forms are of even rank.
\end{proof}

\begin{theorem} $h(3) \ge 6$ and $h(5) \ge 14$.
\end{theorem}
\begin{proof} In general, if $\phi \co \pi_1(M) \to \zz^n$ is an isomorphism, then
$\phi^2 \co H^2(\zz^n) \to H^2(M)$ is injective.

In the case that $n = 3$, all products of two elements  in
$H^2(\zz^3)$ are 0 (since $H^4(\zz^3) = 0$) so  the intersection from on
$H^2(M)$
  vanishes on a rank 3 submodule, implying that this (nonsingular) form must have
rank at least 6.

In the case that $n = 5$, we have the map $H^4(\zz^5) \to H^4(M)
\cong \zz$.  Any such map is given by multiplying with an element $D \in
H^1(\zz^5)$.  After a change of basis, $D$ can be taken to be a multiple of a
generator, say $e_1$.  From this it follows that the intersection form of
$H^2(M)$ vanishes on the 7--dimensional submodule generated by the images of the
set of elements in
$H^2(\zz^5)$,
$\{e_{12}, e_{13},e_{14},e_{15},e_{23},e_{24},e_{25}
\}$ (where $e_{ij} = e_i e_j$).  To see this, observe that the only possible
nontrivial products of two of these are $\pm e_{1234}, \pm e_{1235}$ and
$\pm e_{1245}$, each of which is killed upon multiplying by $e_1$. Since the
nonsingular intersection form on
$H^2(M)$ has a 7--dimensional isotropic subspace, it must be of rank at least 14.
\end{proof}


\section{Algebraic and Geometric 4--Reductions.}

The following algebraic construction will be used repeatedly in constructing our
desired 4--manifolds.

\begin{definition} A {\em 4--reduction of a group $G$ by a 4--tuple of elements}
$[w_1, w_2, w_3, w_4],$
$w_i \in G$, is the quotient of $G$ by normal subgroup generated by the 6
commutators,
$[w_i, w_j], i < j$.  This quotient is denoted $G/[ w_1, w_2, w_3, w_4]$.

 More generally, we say a group $G$ {\em can be 4-reduced to the group
$H$  using  the   4--tuples} $\{[w_{1k},w_{2k}, w_{3k}, w_{4k}]\}_{k=1}^\ell$ if
$H$ is isomorphic to the  quotient of $G$ by the normal subgroup generated by the
$6\ell$ commutators $[w_{ik},w_{jk}], \  i<j, \ k=1,\cdots, \ell$.
\end{definition}

The geometric motivation for this comes from the following theorem.

\begin{theorem}\label{4reduce} If $X$ is a 4--manifold and $\{w_1, w_2, w_3, w_4\}
\subset
\pi_1(X)$, then there is a 4--manifold $X'$ with $\pi_1(X') = \pi_1(X) /[ w_1,
w_2, w_3, w_4]$ and $\beta_2(X') = \beta_2(X) +6.$
\end{theorem}

 Before proving this we make the following simple observation.

\begin{lemma}If a  4--manifold $X'$ is constructed from  a  compact 
4--manifold
$X$ via surgery along a curve
$\alpha$, then $\beta_2(X') = \beta_2(X)$ if $\alpha$ is of infinite order in
$H_1(X)$ and  
$\beta_2(X') =
\beta_2(X) + 2$ otherwise.
\end{lemma}

\begin{proof}Since $X'$ is formed by removing $S^1 \times B^3$ and replacing it
with $B^2 \times S^2$, $\chi(X') = \chi(X) +2$.  If $\alpha $ is of infinite
order, $\beta_1(X') = \beta_1(X) -1$, and similarly $\beta_3(X') = \beta_3(X) -1$
by duality, so $\beta_2(X') = \beta_2(X)$.  On the other hand, if $\alpha$ is of
finite order, $\beta_1$ is unchanged by surgery, by duality $\beta_3$ is unchanged,
so the change in the Euler characteristic must come from an increase in $\beta_2$
by 2.

\end{proof}

\begin{proof}[Proof of Theorem~\ref{4reduce}]

 Form the connected sum $  X \cs T^4$.  This has
increased the second Betti number by six.  Next, perform surgery on four curves to
identify the generators of
$\pi_1(T^4)$ with the elements $w_i$.  This does not change the second Betti
number because the curves being surgered are of infinite order in $H_1(X \cs
T^4)$.  Since the generators of
$\pi_1(T^4)$ commute, the effect of this is that now the four  elements $w_i$
commute. Thus the manifold that results from the surgeries has the stated
properties.
\end{proof}

The main algebraic result concerning 4--reduction, and the key
  to our geometric constructions  via Theorem~\ref{4reduce}, is the following.

\begin{theorem} For $m >2$ and $n >2$, the free product $\zz^m * \zz^n$ can be
4--reduced to $\zz^{m +n}$ using $\frac{mn}{6}$ 4--tuples if $mn$ is divisible by
6.

\end{theorem}

\begin{proof} If the free product $\zz^m * \zz^n$ can be 4--reduced to
$\zz^{m +n}$ using $\frac{mn}{6}$ 4--tuples we will say that the pair
$(m,n)$ is realizable.  Let $\calr$ denote  the set of realizable pairs with 
$m > 2, n>2$.

First we show that $(3,4)$, $(3,6)$, and $(5,6)$ are in $\calr$.
\begin{itemize}

\item Consider first the pair  $(3,4)$. Denote the generators of $\zz^3$    by
$\{x_1, x_2, x_3\}$ and let $\zz^4$ be generated by
$\{y_1, y_2, y_3,y_4\}$.  The two 4--tuples $[x_1, y_1,
x_2y_2, x_3y_3]$ and
$[x_2, x_1y_3, x_3y_2, y_4]$ carry out the desired 4--reduction. For the convenience
of the reader
we provide the details next, but in subsequent examples
similar calculations will be omitted.

  We must show that the subgroup $U$ generated  by these commutator
4--relations contain all 12 commutators $[x_i,y_j]$. It is helpful to recall
 that the set of
elements in a group which commute with a fixed element forms a subgroup.

From the first 4--relation, $[x_1, y_1, x_2y_2, x_3y_3]$, we see, using the
commutators $[x_1, y_1]$, $[x_1 , x_2y_2]$, and $[x_1, x_3y_3]$, that
the commutators
$[x_1,y_1]$, $[x_1,y_2]$, and $[x_1,y_3]$ are in $U$.  The commutators $[
y_1, x_2y_2]$ and
$ [ y_1,  x_3y_3]$ give that the commutators $ [x_2,y_1]$ and $[x_3,y_1]$ are
in
$U$.  The last  commutator, $   [ x_2y_2, x_3y_3]$, we return to momentarily.

From the second 4--relation, $[x_2, x_1y_3, x_3y_2, y_4]$, we have first the
commutator $[x_2, y_3] \in U$.  Then, from  the previous relation ($   [ x_2y_2,
x_3y_3]$) it follows  that      $[x_3, y_2] \in U$. Next, that the
commutators $[x_2, y_2]$ and
$[x_2, y_4]  $ are in $U$ follow immediately.  From the commutator $[  x_1y_3,
x_3y_2]$ we see that $[x_3, y_3] \in U$ (since we already had that $[x_1,
y_2]$ is in $U$).  From
$[x_1y_3, y_4]$ we have $[x_1, y_4] \in U$.  The  commutator $[x_3y_2, y_4]$ gives
the last needed commutator, $[x_3,y_4]$.

\item In the case of the pair $(3,6)$, using similar notation, the following three
4--tuples
$[x_1, y_1, x_2y_2, x_3y_3], [x_2, x_1y_3, x_3y_4, x_2y_5], [y_6, x_1y_2, x_2x_3,
x_3y_4]$ reduce $\zz^3 * \zz^6$.

\item Finally, for     $(5,6)$, using similar notation, the
following five 4--tuples suffice:
$[x_1,y_1,x_2y_2,x_3y_3]$,
$[x_2,y_3,x_3y_4,x_4y_5]$,
$[x_5, y_6,  x_4y_3,x_1y_2y_5 ]$,
$[x_3,y_5,x_5y_6,x_1y_4] $,
$[x_1y_1, x_2y_4, x_4y_6, x_5y_2]$.

\end{itemize}
For the general case of $(m,n)$, assume first that $m$ is divisible by 6.  Using
the realization of
$(3,4)$ we can realize $(6,4)$ and have already realized
$(6,3)$.  (Separate the six generators into two groups of three and make each set
commute with the other four using the construction used for $(3,4)$.)  Combining
these, we can realize
$(6k,3)$ and
$(6k,4)$ for any
$k$.  Now, combining these we can realize $(6k, 3a + 4b)$ for any $a$ and
$b$.  But all integers greater than 2, other than 5, can be written as
$3a + 4b$ for some $a$ and $b$.

In the case that neither $m$ nor $n$ are divisible by 6, we can assume 3 divides
$m$ and we want to realize $(3k, n)$.  Notice that $n$ must be even.  Since we
can realize
$(3,4)$ and $(3,6)$, we can realize $(3k,4)$ and
$(3k,6)$ for all $k$.  Thus we can realize $(3k, 4a + 6b)$ for all $a$ and $b$,
but $4a + 6b$ realizes all even integers greater than 3.
\end{proof}


\section{Basic Realizing Examples}

We begin with the exceptional cases of $n=3$ and $n=5$ and then move on to a set
of basic examples for which $h(n)  = C(n,2)+\epsilon_n$.    In the next section we
note that these basic examples can be used to construct the necessary examples
for the proof of Theorem~\ref{mainthm}.

\begin{itemize}
\item{\bf $n=3$}:  Start with the 4--torus, $T^4$, with $\beta_2(T^4) = 6$.
Surgery on a single curve representing a generator of $\pi_1(T^4)$ results in a
manifold
$M$ with
$\pi_1(M) = \zz^3$ and $H_2(M) =
\zz^3$.

\item {\bf $n=5$}:  Begin with $X= F_2 \times F_4$ with $\beta_2(X) = 10$
and $\pi_1$ generated by
$\{x_1,x_2\}$ and $\{y_1, y_2, y_3, y_4\}$.    Perform a surgery to identify
$y_3$ and $y_4$, so that the group is generated by $\{x_1, x_2, y_1, y_2, y_3
\}$.  Notice that $y_1$ and $y_2$ commute, as follows from the original surface
commutator relationship
$[y_1, y_2][y_3, y_4] =1$. This surgery, since it is along an element of infinite
order in
$H_1$, does not change $H_2(X)$. Hence, it only remains to arrange that the pairs
of elements $\{y_1, y_3\}$ and $\{y_2, y_3\}$ commute.  Performing
surgery on a (rationally) null homologous curve raises $\beta_2$ by two,
so performing surgeries to kill these two commutators raises the rank of
$H_2(X)$ by 4, and   the resulting 4--manifold $M$ has
$H_2(M) = \zz^{14}$ as desired.

\end{itemize}

We will say that  the integer {\em $n$ is realizable} if there is a closed
oriented 4--manifold $M_n$ with
$\pi_1(M_n) =
\zz^n$ and $\beta_2(M_n) = C(n,2)+\epsilon_n$.  Let $\cals$ be the set of
realizable integers.
We now show that ${0,1,2,4,6,7,8,9,11,12} \subset \cals$ by describing the
construction of  realizing 4--manifolds $M_n$ for each of these $n$.

\begin{itemize}
\item{\bf $n=0$}:  $M_0 = S^4$.

\item{\bf $n=1$}:  $M_1 = S^1 \times S^3$.

\item{\bf $n=2$}:  $M_2 = F_2 \times S^2$. Notice that $C(2,2)=1$, so that
$\beta_2(M_2)=2=C(2,2)+\epsilon_2$.

\item{\bf $n = 4$}:  $M_4 = T^4$.

\item{\bf $n=6$}:  Build $M_6$ as follows. Let $X = F_2 \times F_4$ with
$\pi_1$ generated by  the 6 elements $\{x_1, x_2 \}$, $\{y_1, y_2, y_3,
y_4\}$.   Note that
$\beta_2(X)= 10$.  Apply Theorem~\ref{4reduce} to
perform the 4--reduction  $[y_1, y_2, y_3, y_4]$ and arrive at the 4--manifold
$M_6$ with $\pi_1(M_6) = \zz^6$ and $\beta_2(M_6) = 16 = C(6,2)+\epsilon_6$ as
desired.

\item{\bf $n=7$}: Begin with   $X = F_2 \times F_4 \cs T^4 $ so  $\beta_2(X) =
16$. Let the generators of $\pi_1$ be
$\{x_1, x_2\}$, $\{y_1, y_2, y_3, y_4\}$, and
$\{z_1, z_2, z_3, z_4\}$. Perform surgeries giving
$y_1=z_2$, $y_2=z_3$, $y_4=z_4$. Now use Theorem~\ref{4reduce} to perform the
4--tuple reduction
$[z_1, x_1y_2, x_2y_1, y_3]$.  (In checking that this abelianizes the group, use
the fact the
$[y_1,y_2]=1$ if and only if
$[y_3,y_4] =1$.)  The resulting 4--manifold $M_7$ has $\pi_1(M_7) = \zz^7$ and
$\beta_2(M_7) = 22 = C(7,2) + \epsilon_7$.

\item{\bf $n = 8$}: Take  $X =  F_4 \times F_4  \cs  F_2 \times F_4$, with
$\beta_2(X)= 28$ and  $\pi_1$ generated by $\{x_1, x_2, x_3, x_4\}$, $\{y_1, y_2,
y_3, y_4\}$, 
$\{z_1, z_2\}$ and $ 
 \{w_1, w_2, w_3,w_4\}$, respectively. Now perform surgeries to introduce the
following relations:

$z_1 = x_1 y_1$,

$z_2 = x_2 y_2$,

$w_1 = x_3 y_2$,

$w_2 = x_4$,

$w_3 = x_1 y_3$,

$w_4 = x_2 y_4$.

\noindent Since $[w_1,w_2][w_3,w_4] = 1$ we have that $[x_3y_2,
x_4][x_1y_3,x_2y_4] =1$.  Since the $x_i$ commute with the $y_i$, this implies
that $[x_3,x_4][x_1,x_2][y_3,y_4] =1$.  From the surface relation for the $x_i$ 
$[x_1,x_2][x_3,x_4] =1$,   it then follows that
$[y_3,y_4] = 1$.  From this and the surface relation for the $y_i$ it follows that
$[y_1, y_2]=1$. Now, the fact that $[z_1, z_2] =1$ gives that $[x_1,x_2] =1$, which
implies that
$[x_3,x_4] =1$ too. The remaining needed four commutator relations between the
$x_i$  and the four commutator relations between the $y_i$ gives  a total of 8
needed commutator relations.  These are gotten by first considering the relations
$[z_1, w_i]$ and then the relations
$[z_2, w_i]$. The resulting $M_8 $ has $\pi_1(M_8) = \zz^8$ and  $\beta_2(M_8) =
28 = C(8,2) +
\epsilon_8$.

\item{\bf $n=9$}: Start with $X = T^4 \cs T^4 \cs T^4$, with $\beta_2(X) =
18$ and $\pi_1$ generated by
$\{x_1, x_2, x_3, x_4\}$, $\{y_1, y_2, y_3, y_4\}$, and
$\{z_1, z_2, z_3, z_4\}$.  Perform three surgeries to give the identifications:
$z_2 = x_3y_3$, $z_3 = x_4$, and $z_4 = y_4$.   Notice that since the 4--tuple
relation $[z_1, z_2, z_3, z_4]$ held in the original group, we now have the
4--tuple relation
$[z_1 , x_3y_3, x_4, y_4]$.

Use Theorem~\ref{4reduce} to add the following three more 4--tuple relations,
raising
$\beta_2$ to 36:

$[ z_1, y_1, x_3, x_2y_2     ],$

$[ z_1 x_3, x_1, y_2, x_2y_3   ],$

$[  x_2y_4, x_4y_3, x_2y_1, x_1y_2y_4    ].$

\noindent The resulting manifold $M_9$ has $\pi_1(M_9) = \zz^9$ and $\beta_2(M_9)
= 36 = C(9,2) +
\epsilon_9$.

\item{\bf $n = 11$}: Start with $X= F_4 \times F_6 \cs T^4 $, with $\beta_2(X)
= 32$ and
$\pi_1$ generated by
$\{x_1, x_2, x_3, x_4\}$, $\{y_1, y_2, y_3, y_4, y_5, y_6\}$, and $\{z_1, z_2,
z_3, z_4\}$.  Perform surgery to get the following identifications:

 $y_4 = z_2$,

 $y_5 = z_3$,

$y_6 = z_4$.

\noindent This leaves a generating set with eleven elements:  $$\{ x_1, x_2, x_3,
x_4, y_1, y_2, y_3, y_4, y_5, y_6, z_1\}.$$  (Notice that $[x_1, x_2][x_3,x_4] =1$
and
$[y_1,y_2][y_3,y_4] = 1$, since $y_5$ and $y_6$ now commute.) Apply
Theorem~\ref{4reduce} four times to perform the following 4--reductions:

$[x_1, x_2, x_3, z_1],$

$[y_1, y_2, y_3, z_1],$

$[x_4y_1 , y_5, y_3y_6, y_2z_1],$

$[x_1y_4, x_4 y_6, x_2 y_2 z_1, y_1 y_3].$

\noindent The resulting manifold $M_{11}$ has $\pi_1(M_{11}) = \zz^{11}$ and
$\beta_2(M_{11}) = 32 + 24 = 56 = C(11,2) + \epsilon_{11}$.

\item{\bf $n = 12$}:  Start with $X = F_4 \times F_4 \cs T^4 $ with $\beta_2(X) =
24$ and $\pi_1$ generated by $\{x_1, x_2, x_3, x_4\}$, $\{y_1, y_2, y_3, y_4\}$, and
$\{z_1, z_2, z_3, z_4\}$ as before. Now apply
Theorem~\ref{4reduce} to add   seven 4--tuple relations:

$[x_1, x_2, x_3, z_1],$

$[y_1, y_2, y_3, z_1],$

$[x_1, x_4, z_2, y_1z_3],$

$[y_1, y_4, z_3, x_3z_4],$

$[x_2y_2, x_4, y_4z_1, z_1z_4],$

$[x_2z_3, z_4, y_3, x_3z_2],$

$[ x_3z_3, y_4z_2,x_1y_2, z_4x_2z_2].$

\noindent The resulting $M_{12} $ has $\beta_2(M_{12}) = 24 + 42 = 66 = C(12,2) +
\epsilon_{12}$ and $\pi_1(M_{12}) = \zz^{12}$ as desired.
\end{itemize}


\section{Constructing more examples: the proof of Theorem~\ref{mainthm}}

\begin{theorem}\label{real} If $m \in \cals$, $n \in \cals$, $(m,n) \in
\calr$,  and if one of
$m, \  m-1, \  n, $ or $\ n-1$ is congruent to $0$ modulo $4$, then $m+n
\in \cals$.
\end{theorem}

\begin{proof} The stated mod 4 condition together with the fact that
$mn\equiv 0 \pmod{6}$ assures that $C(n+m,2)+\epsilon_{n+m}=(C(n,2)+\epsilon_n)
+(C(m,2)+\epsilon_m)+mn$. Thus, one can build the desired $M_{m+n}$ by performing
$mn$ surgeries on
$M_m \cs M_n \cs \frac{mn}{6}T^4$; that is, by performing $\frac{mn}{6}$
4--reductions as in Theorem~\ref{4reduce}.

\end{proof}

We have that $\{0,1,2,4,6,7,8,9,11,12\} \subset \cals$. Furthermore, all pairs
$(n,m)$ with
$mn \equiv 0 \mod 6$ and $m \ge 3$ and $n \ge 3$ are in $\calr$.

Using Theorem~\ref{real} and the pair $(4,6) \in \calr$ gives $10 \in
\cals$.  Similarly, using the pair $(4,9) \in \calr$ gives $13 \in
\cals$.  The pair $(6,8) \in \calr$ shows that $14 \in \cals$.  The pair
$(6,9) \in \calr$ shows that $15 \in \cals$.   The pair $(4,12) \in
\calr$ shows that $16
\in \cals$.   The pair $(8,9) \in \calr$ shows that $17
\in \cals$.

Next the pairs $(12,n)$ for $n = 6, \ldots , 17$ show that $\{18, \ldots , 29\}
\subset \cals$.  Then the pairs  $(12,n)$ for $n = 18, \ldots , 29$ show that
$\{30,
\ldots , 41\}
\subset \cals$. Repeating inductively in this way shows that $n \in
\cals$ for all $n \ge 6$, as desired.


\section{Finitely generated abelian groups}\label{finitesec}

We next use the manifolds constructed in the previous sections as building blocks
  to prove Theorem \ref{genral}.  Let
$$G=\zz/d_1\oplus \cdots \zz/d_k\oplus \zz^n$$ where $d_i|d_{i+1}$.
Fix a prime $p$ that divides $d_1$. Then
\begin{equation*}
\rk_{\zz/p}H_1(G;\zz/p)=k+n \text{ and } \rk_{\zz/p}H_2(G;\zz/p)= C(k+n,2)+k.
\end{equation*} If $X$ is a closed, oriented 4--manifold with $\pi_1(X)\cong G$,
we have   $H_1(X;\zz/p)\cong H_1(G;\zz/p)$ and $H_2(X;\zz/p)$ surjects
to
$ H_2(G;\zz/p)$.  This gives the lower bound on $q(G)$
\begin{equation}\label{torsionbound} q(G)\ge  2-2(k+n)+ C(k+n,2)+k=1-n+
C(n+k-1,2).\end{equation}

To construct upper bounds, consider the following constructions
 of closed 4--manifolds $X$ with
$\pi_1(X)\cong G$. Let $L(d)$ denote a 3-dimensional lens space with
$\pi_1(L(d))=\zz/d$, and let $B_n$ denote a closed 4--manifold with
$\pi_1(B_n)=\zz^n$ and $\chi(B_n)=q(\zz^n)$.

\begin{itemize}

  \item  Let   $X$ be the manifold obtained by starting with $B_{k+n}$ and doing
surgeries on the
  first $k$ generators of $\pi_1(B_{k+n})=\zz^{k+n}$ in such a way as to kill
$d_1$ times the first generator, $d_2$ times the second, and so forth. Then
$\pi_1(X)=\zz/d_1\oplus\cdots \oplus\zz/d_k \oplus \zz^n$ and
$\rk_\qq(H_2(X;\qq))= \rk_\qq(H_2(B_{n+k};\qq))$.
  Therefore
  $$ \chi(X)=2-2n +\rk_\qq(H_2(B_{n+k};\qq)). $$
  Thus, if $n+k\ne 3,5$,   simplifying yields
\[
  \chi(X)=1 -n +C(n+k-1,2) + k + \epsilon_{n+k}
\] and hence
 \begin{equation}\label{torsionbound2}
 0\leq q(G)-\big(1-n+ C(n+k-1,2)\big) \leq  k +\epsilon_{n+k}.\end{equation}
\item  Suppose that $n\ge 1$.  Start with
$$Y=\big((L(d_1)\# L(d_2)\#\cdots \# L(d_k))\times S^1\big)\# B_{k+n-1}.$$ Then
perform $k$ surgeries which identify the $k$ generators of the connected sum of
lens spaces with the first $k$ generators of $\pi_1(B_{k+n-1})$. These surgeries
do not change the  rank of
$H_2(Y;\qq)$. Finally perform $n-1$ surgeries along circles representing the
commutator of the $S^1$ factor and the last $n-1$ generators of
$\pi_1(B_{k+n-1})$.  Each of these surgeries  increase the rank of the second
rational homology by 2 since the commutators are nullhomologous.   This produces
a 4--manifold $X$ with $\pi_1(X)\cong G$ and
$\rk_{\qq}H_2(X;\qq)=2(n-1) + \rk_{\qq}H_2(B_{k+n-1};\qq)$. Hence
\[
\chi(X)= \rk_{\qq}H_2(B_{k+n-1};\qq).\] Thus, if $n+k-1\ne 3, 5$,
$
\chi(X)=  C(n+k-1,2)+\epsilon_{n+k-1}.
$

Referring to Equation (\ref{torsionbound}),  this gives the upper bound
    \begin{equation}\label{torsionbound1}
 0\leq q(G)-\big(1-n+ C(n+k-1,2)\big) \leq n-1
+\epsilon_{n+k-1} \end{equation}
  (for $n,k$ satisfying  $n+k-1 \ne 3,5$ and  $n\ge 1$).

 \item Consider now the case $n=0$. Start with the 4--manifold obtained from
 $\big( (L(d_1)\#\cdots \#L(d_{k-1}))\times S^1\big) \# B_{k-1}$ by
performing surgery to identifying the generators of $ \pi_1(B_{k-1})$ with the
 generators for the lens spaces. This yields a closed 4--manifold $Y$ with
 $\pi_1(Y)\cong \zz/d_1\oplus\cdots\oplus \zz/d_{k-1}\oplus \zz$. Surgering $d_k$
times the last generator  gives a  closed 4--manifold $X$ with
 $\pi_1(X)\cong \zz/d_1\oplus\cdots\oplus \zz/d_{k}$ and
$\rk_\qq(H_2(X;\qq))=\rk_\qq(H_2(B_{k-1};\qq))$.  Therefore
 \[\chi(X)=2+\rk_\qq(H_2(B_{k-1};\qq)).\]
 Thus  when $k-1\ne 3,5$,
$\chi(X)=2+C(k-1,2) +\epsilon_{k-1},$
 and so  for $G=\zz/d_1\oplus\cdots\oplus \zz/d_{k}$,
 \begin{equation}\label{torsionbound3} 0\leq q(G)-\big(1+ C(k-1,2)\big) \leq  1
+\epsilon_{k-1}.
\end{equation}

\item In two cases the rational homology gives better lower bounds than the $\zz/p$
homology. For
 $G=\zz/d_1\oplus \cdots\oplus \zz/d_k$,
 $\rk_\qq(H_1(G;\qq))=0$ and $\rk_\qq(H_2(G;\qq))=0 $ and hence $q(G)\ge 2$. Combined with
 Equation (\ref{torsionbound3})
 this shows that $q(\zz/d)= 2$ and
$q(\zz/d_1\oplus\zz/d_2)= 2$.
 \end{itemize}


The estimates (\ref{torsionbound}), (\ref{torsionbound1}),
(\ref{torsionbound2}),  and (\ref{torsionbound3}) and the discussion of the  previous paragraph combine to give a proof of Theorem
  \ref{genral}.

\section{Remarks and questions}

A variant of $q(G)$ is obtained by defining $p(G)$ to be the smallest value of $\chi(X)- \sigma(X) $
for all closed oriented  4--manifolds $X$ with $\pi_1(X)=G$. Here $\sigma(X)$ denotes the signature of $X$.
Notice that all the examples we constructed for $G$ abelian have signature zero. We conjecture that
$p(\zz^n)=2-2n+C(n,2)$. This guess is motivated by a slight amount of redundancy which occurs in the constructions given above when $\epsilon_n=1$.  For example, in the case of $n=6$, the construction we gave starts with $F_2\times F_4$ and uses one 4--reduction to abelianize the fundamental group of $F_4$. In particular, the  surface relation $[y_1,y_2][y_3,y_4]=1$ shows that one of the 6 relations coming from the 4--reduction is unnecessary. This extra bit of flexibility 
may perhaps be used to twist the geometric construction slightly to introduce some signature.

 An interesting question is whether the invariant $q$ depends on the category of manifold
 or choice of geometric structure.
 For example, one might consider the infimum of $\chi(X)$ over smooth 4--manifolds or topological 4--manifolds
or even 4--dimensional  Poincar\'e complexes with $\pi_1(X)=G$. The manifolds we constructed for $G=\zz^n$ are smooth and
the lower bounds are homotopy invariants, so that for $G=\zz^n$ the value of $q$ is independent of the category.


\newcommand{\etalchar}[1]{$^{#1}$}

\end{document}